\pdfoutput=1    

\documentclass[11pt,a4paper]{article}

\usepackage{amsmath}
\usepackage{amssymb,verbatim}
\usepackage{amsfonts}
\usepackage{amsthm}
\usepackage{pifont}
\usepackage{graphicx}
\usepackage{datetime}
\usepackage{scrdate}
\usepackage{enumitem}
\usepackage{geometry}
\usepackage{float}
\usepackage{bbm} 

\newtheorem{theorem}{Theorem}

\theoremstyle{definition}
\newtheorem{example}{Example}

\newcommand{\ds}{\displaystyle}
\newcommand{\beq}{\begin{equation*}}
\newcommand{\eeq}{\end{equation*}}
\newcommand{\beqn}{\begin{equation}}
\newcommand{\eeqn}{\end{equation}}
\newcommand{\R}{\mathbb{R}}
\newcommand{\E}{\mathbb{E}}

\newcommand{\dd}{\mathrm{d}}
\newcommand{\Ol}{\varOmega}
\newcommand{\ph}{\varphi}
\newcommand{\thet}{\vartheta}
\newcommand{\rh}{\varrho}
\newcommand{\eps}{\varepsilon}
\newcommand{\tn}{\textnormal}

\begin{document}

\title{The mean width of the oloid and\\ integral geometric applications of it}
\author{Uwe B\"asel}
\date{}
\maketitle

\begin{abstract}
\noindent The oloid is the convex hull of two circles with equal radius in perpendicular planes so that the center of each circle lies on the other circle. We calculate the mean width of the oloid in two ways, first via the integral of mean curvature, and then directly. Using this result, the surface area and the volume of the parallel body are obtained. Furthermore, we derive the expectations of the mean width, the surface area and the volume of the intersections of a fixed oloid and a moving ball, as well as of a fixed and a moving oloid.\\[0.2cm]
\end{abstract}

\section{Introduction}

\begin{figure}[h]
\begin{minipage}[c]{0.5\textwidth}
\begin{center}
\includegraphics[scale=0.301]{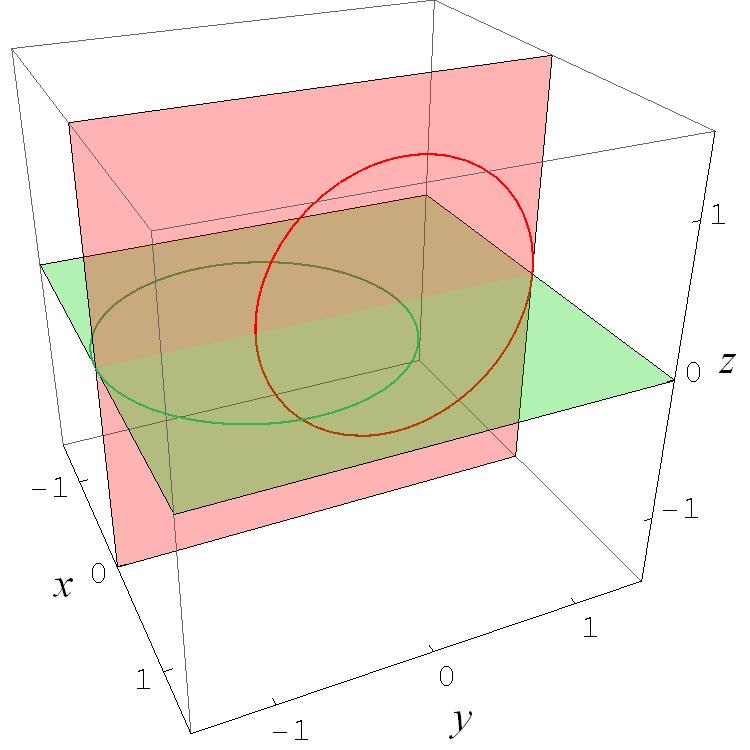}
\caption{\label{Fig:circles} The circles $k_A$ and $k_B$ with $r=1$}
\end{center}
\end{minipage}
\hspace{0.2cm}
\begin{minipage}[c]{0.5\textwidth}
\begin{center}
\includegraphics[scale=0.30]{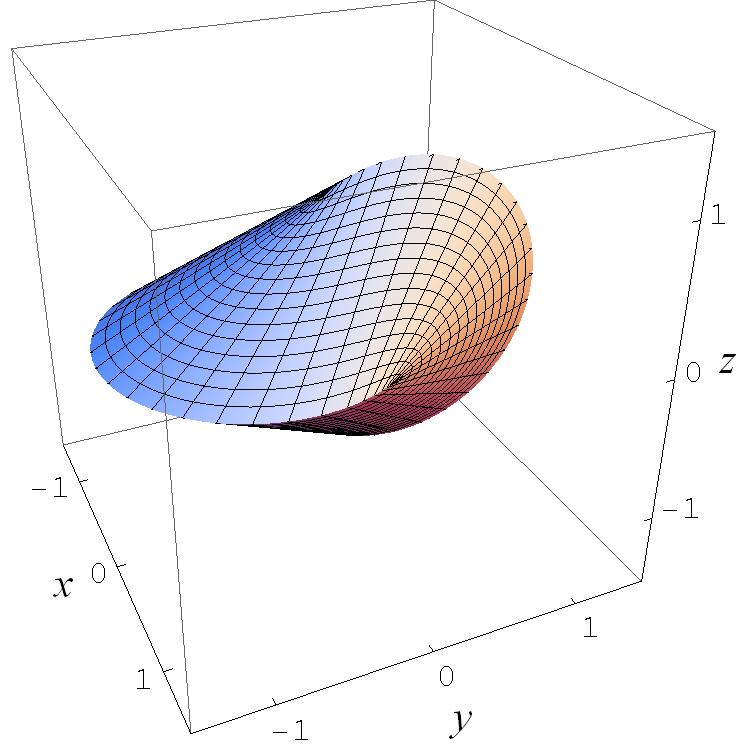}
\caption{\label{Fig:oloid} The oloid $\Ol_1$}
\end{center}
\end{minipage}
\end{figure}

The oloid $\Ol_r$ is the convex hull of two circles $k_A$, $k_B$ with equal radius $r$ in perpendicular planes so that the center of each circle lies on the other circle (see Figures \ref{Fig:circles} and \ref{Fig:oloid}). Dirnb\"ock \& Stachel \cite[p.\ 117]{Dirnboeck_Stachel} calculated the surface area and the volume of the oloid (see also \cite{Wikipedia}, \cite{Oloide}, and Equations \eqref{Eq:S(O_1)-1}, \eqref{Eq:S(O_1)-2}, \eqref{Eq:V(O_1)-1} and \eqref{Eq:V(O_1)-2} of the present paper). The surface $\partial\Ol_r$ is part of a developable surface \cite{Dirnboeck_Stachel}, \cite{Baesel_Dirnboeck}, \cite{Wikipedia}, \cite{Oloide}.

Finch \cite{Finch} calculated surface areas, volumes and mean widths of the convex hulls of three different configurations of two orthogonal disks with equal radius. The mean width $\bar{b}$ of every convex hull is determined twice: 1) using the integral $M$ of the mean curvature and the relation $\bar{b} = M/(2\pi)$, 2) calculating $\bar{b}$ directly.\\[0.2cm]   
\hspace*{0.4cm} According to \cite[pp.\ 105-106]{Dirnboeck_Stachel}, the circles with $r=1$ can be defined by
\beqn \label{circles}
\begin{aligned}
  k_A
:= {} & \left\{(x,y,z)\in\R^3 \colon x^2+(y+1/2)^2=1 \wedge z=0\right\},\\
  k_B
:= {} & \left\{(x,y,z)\in\R^3 \colon (y-1/2)^2+z^2=1 \wedge x=0\right\}
\end{aligned}
\eeqn
(see Fig.\ \ref{circles}). A parametrization of the surface $\partial\Ol_1$ \cite[p.\ 165, Eq.\ (2)]{Baesel_Dirnboeck} is
\beqn \label{Eq:parametrization_1}
  (x,y,z)^\tn{T}
= \vec{\omega}(m,t) 
= (\omega_1(m,t),\omega_2(m,t),\pm\omega_3(m,t))^\tn{T}\,,\quad
  0\le m\le 1\,,\quad -\frac{2\pi}{3}\le t\le\frac{2\pi}{3}\,,
\eeqn
with
\beqn \label{Eq:parametrization_2}
\left.
\begin{aligned}
  \omega_1(m,t)
= {} & (1-m)\sin t\,,\\[0.1cm]
  \omega_2(m,t)
= {} & \frac{2(m-1)\cos^2 t+(2m-3)\cos t+2m-1}{2(1+\cos t)}\,,\\
  \omega_3(m,t)
= {} & \frac{m\sqrt{1+2\cos t}}{1+\cos t}\,.
\end{aligned}
\;\;\right\}
\eeqn
\hspace*{0.4cm} To the authors knowledge, the mean width of the oloid is not aready known. In Section~\ref{S:Mean_curvature} we calculate the mean width of $\Ol_r$ using the integral of mean curvature, and in Section \ref{S:Mean_width} we calculate it directly. With the help of this result we derive the volume, the surface area and the mean width of the parallel body of $\Ol_r$ in Section \ref{S:Parallel_body}. Using the principal kinematic formula of integral geometry, the expectations of the mean width, the surface area and the volume of the intersections of a fixed oloid and a moving ball, as well as of a fixed and a moving oloid are calculated in Section \ref{S:Intersections}.

\section{Preliminaries}

In the following, we work in the real vector space $\R^3$ with its standard scalar product $\big\langle\vec{a},\vec{b}\,\big\rangle = \vec{a}\,\cdot\,\vec{b}$ and its vector product $\vec{a}\times\vec{b}$ for two vectors $\vec{a}=(a_1,a_2,a_3)^\tn{T}$, $\vec{b}=(b_1,b_2,b_3)^\tn{T}$. We denote the partial derivatives
\beq
  \frac{\partial\vec{\omega}}{\partial m}\,,\quad
  \frac{\partial\vec{\omega}}{\partial t}\eeq
 of $\vec{\omega}=\vec{\omega}(m,t)$ (see \eqref{Eq:parametrization_1}) by $\vec{\omega}_m$, $\vec{\omega}_t$, and so on.\\[0.2cm]  
\hspace*{0.4cm} Using \eqref{Eq:parametrization_2}, for the coefficients $g_{11}=E$, $g_{12}=F=g_{21}$, $g_{22}=G$ of the first fundamental form (see e.\,g. \cite[pp.\ 87-88]{Kreyszig}, translation: p.\ 68) we find
\beqn \label{Eq:first_fundamental_form}
\left.
\begin{aligned}
  g_{11}
= {} & \langle\vec{\omega}_m,\vec{\omega}_m\rangle
= 3\,,\qquad
g_{12}
= \langle\vec{\omega}_m,\vec{\omega}_t\rangle
= \tan(t/2)\,,\\[0.1cm]
  g_{22}
= {} & \langle\vec{\omega}_t,\vec{\omega}_t\rangle
= \frac{2(3m^2-4m+1)\cos^2 t-(4m-3)\cos t+1}{(1+\cos t)(1+2\cos t)}\,,\\[0.1cm]
  g
= {} & \det(g_{jk})
= \left|\begin{array}{ll}
	g_{11} & g_{12}\\
	g_{21} & g_{22}
  \end{array}\right|
= g_{11}\,g_{22}-g_{12}^2
= \frac{2[(3m-2)\cos t-1]^2}{(1+\cos t)(1+2\cos t)}\,.
\end{aligned}
\;\;\right\}
\eeqn
Now, we able to calculate the surface area of the oloid $\Ol_1$:
\begin{align} \label{Eq:S(O_1)-1}
  S(\Ol_1)
= {} & \int_{\partial\Ol_1}\dd S
= \int_{\partial\Ol_1}\dd S(m,t)
= 2\int_{t=-2\pi/3}^{2\pi/3}\,\int_{m=0}^1\sqrt{g(m,t)}\:\dd m\,\dd t\nonumber\displaybreak[0]\\[0.2cm]
= {} & 4\int_{t=0}^{2\pi/3}\int_{m=0}^1\sqrt{g(m,t)}\:\dd m\,\dd t
= 4\int_{t=0}^{2\pi/3}\int_{m=0}^1\frac{\sqrt{2}\left[(2-3m)\cos t+1\right]}{\sqrt{(1+\cos t)(1+2\cos t)}}\,\dd m\,\dd t\nonumber\displaybreak[0]\\
= {} & 2\,\sqrt{2}\int_0^{2\pi/3}\frac{2+\cos t}{\sqrt{(1+\cos t)(1+2\cos t)}}\,\dd t\,.
\end{align}
{\em Mathematica} evaluates this integral to
\beqn \label{Eq:S(O_1)-2}
  S(\Ol_1)
= 2\,\sqrt{2}\cdot\sqrt{2}\,\pi
= 4\pi\,.
\eeqn
\hspace*{0.4cm} Now, we calculate the volume of $\Ol_1$, and start with
\begin{align*}
  V(\Ol_1)
= {} & 2\iint z\,\dd x\,\dd y
= 2\int_{t=-2\pi/3}^{2\pi/3}\,\int_{m=0}^1\omega_3(m,t)\left|\frac{\partial(\omega_1(m,t),\,\omega_2(m,t))}{\partial(m,t)}\right|\dd m\,\dd t\\
= {} & 4\int_{t=0}^{2\pi/3}\,\int_{m=0}^1\omega_3(m,t)\left|\frac{\partial(\omega_1(m,t),\,\omega_2(m,t))}{\partial(m,t)}\right|\dd m\,\dd t\,.
\end{align*}
From \eqref{Eq:parametrization_2} it follows that
\beq
  \frac{\partial(\omega_1(m,t),\,\omega_2(m,t)}{\partial(m,t)}
= \left|
  \begin{array}{ll}
	\ds{\frac{\partial\omega_1(m,t)}{\partial m}} & 
	\ds{\frac{\partial\omega_1(m,t)}{\partial t}}\\[0.3cm]
	\ds{\frac{\partial\omega_2(m,t)}{\partial m}} &
	\ds{\frac{\partial\omega_2(m,t)}{\partial t}}
  \end{array}
  \right|
= -\frac{1+(2-3m)\cos t}{1+\cos t}\,.
\eeq
So we have   
\begin{align} \label{Eq:V(O_1)-1}
  V(\Ol_1)
= {} & 4\int_{t=0}^{2\pi/3}\int_{m=0}^1\frac{m\,\sqrt{1+2\cos t}}{1+\cos t}\,\frac{1+(2-3m)\cos t}{1+\cos t}\,\dd m\,\dd t\nonumber\\[0.1cm]
= {} & 2\int_0^{2\pi/3}\frac{\sqrt{1+2\cos t}}{(1+\cos t)^2}\,\dd t\,.
\end{align}
{\em Mathematica} finds
\beqn \label{Eq:V(O_1)-2}
  V(\Ol_1)
= \frac{2}{3}\left[K\big(\sqrt{3}\big/2\big)+2E\big(\sqrt{3}\big/2\big)\right],
\eeqn
where
\beqn \label{Eq:first_kind}
  K(k)
= F(\pi/2,k)
= \int_0^{\pi/2}\frac{\dd x}{\sqrt{1-k^2\sin^2 x}}\,,
\eeqn
is the complete elliptic integral of the first kind, and
\beqn \label{Eq:second_kind}
  E(k)
= E(\pi/2,k)
= \int_0^{\pi/2}\sqrt{1-k^2\sin^2 x}\,\dd x  
\eeqn
is the complete elliptic integral of the second kind. A numerical integration integration of \eqref{Eq:V(O_1)-1} and evaluation of \eqref{Eq:V(O_1)-2} with {\em Mathematica} yields the decimal expansion
\beq
   V(\Ol_1)
\approx  3.05241846842437485669720053193  
\eeq
(see also \cite{Alcover}). 

\section{The integral of mean curvature} \label{S:Mean_curvature}

The surface $\partial\Ol_1$ of the oloid $\Ol_1$ is piecewise continuously differentiable. We denote by $H$ the mean curvature in one point of $\partial\Ol_1$. The circles $k_A$ and $k_B$ (see \eqref{circles}) produce two edges $\eps_1$ and $\eps_2$, respectively, that are smooth curves. Let $\alpha=\alpha(t)$ denote the angle between the two unit normal vectors in every point of $\eps_1$. Applying the general formula for the integral $M$ of the mean curvature (see \cite[pp.\ 76-77, Eqs.\ (3.5), (3.7)]{Voss}; cp.\ the formula for the mean width in \cite[p.\ 3]{Finch}) to $\Ol_1$ gives
\beqn \label{Eq:M(O_1)}
  M(\Ol_1)
= \int_{\partial\Ol_1}H\,\dd S 
+ \frac{1}{2}\,\sum_{j=1}^2\int_{\eps_j}\alpha\,\dd s
= \int_{\partial\Ol_1}H(m,t)\:\dd S(m,t) 
+ \int_{\eps_1}\alpha(t)\,\dd t
\eeqn
For the unit normal vector one finds
\beqn \label{n-vec}
  \vec{n} 
= (n_1,n_2,n_3)^\tn{T}
= \frac{\vec{\omega}_m\times\vec{\omega}_t}
       {\lvert\vec{\omega}_m\times\vec{\omega}_t\rvert}
= \left(\sin(t/2)\,,\:-\frac{\cos t}{2\cos(t/2)}\,,\:
	\frac{\sqrt{1+2\cos t}}{2\cos(t/2)}\,\right)^\tn{\!\!T}.	 
\eeqn
Since $\partial\Ol_1$ is part of a developable surface and the line segments $t=\tn{const}$, $0\le m\le 1$ are part of the generators of $\partial\Ol_1$ (see \cite{Dirnboeck_Stachel}, \cite{Baesel_Dirnboeck}), it is not surprising that $\vec{n}$ does not depend on $m$. The mean curvature in a point of a surface is defined by
\beq
  H 
= \frac{1}{2}\,(\kappa_1+\kappa_2)
= \frac{1}{2g}\,(g_{11}b_{22}-2g_{12}b_{12}+g_{22}b_{11})\,,
\eeq
where $\kappa_1$, $\kappa_2$ are the principal curvatures, $b_{ik}$ are the coefficients of the second fundamental form (see e.\,g. \cite[p.\ 99]{Kreyszig}, translation: p.\ 79), and $g_{ik}$, $g$ are given by \eqref{Eq:first_fundamental_form}. In our case we have  
\beq
  H(m,t)\,\dd S(m,t)
= H(m,t)\,\sqrt{g(m,t)}\,\dd m\,\dd t
= \frac{1}{2\,\sqrt{g}}\,(g_{11}b_{22}-2g_{12}b_{12}+g_{22}b_{11})\,\dd m\,\dd t 
\eeq
and
\begin{align*}
  b_{11} 
= {} & L
= \langle\vec{\omega}_{mm},\vec{n}\rangle
= 0\,,\qquad
  b_{12} 
= M
= \langle\vec{\omega}_{mt},\vec{n}\rangle
= 0\,,\\[0.1cm]
  b_{22} 
= {} & N
= \langle\vec{\omega}_{tt},\vec{n}\rangle
= \frac{(3m-2)\cos t-1}{\sqrt{2}\,(1+2\cos t)\,\sqrt{1+\cos t}}\,.
\end{align*}
It follows that
\begin{align*}
  H(m,t)\,\dd S(m,t)
= {} & \frac{g_{11}(m,t)b_{22}(m,t)}{2\,\sqrt{g(m,t)}}\:\dd m\,\dd t\\
= {} & \frac{3}{2}
  \cdot\frac{(3m-2)\cos t-1}{\sqrt{2}\,(1+2\cos t)\,\sqrt{1+\cos t}}
  \cdot\frac{\sqrt{1+\cos t}\,\sqrt{1+2\cos t}}{\sqrt{2}\,[(3m-2)\cos t-1]}\,\dd m\,\dd t\\
= {} & \frac{3}{4\,\sqrt{1+2\cos t}}\,\dd m\,\dd t	 
\end{align*}
and
\begin{align*}
  \int_{\partial\Ol_1}H\,\dd S
= {} & 2\int_{m=0}^1\int_{t=-2\pi/3}^{2\pi/3}H(m,t)\,\dd S(m,t)
= \frac{3}{2}\int_0^1\dd m\:\int_{-2\pi/3}^{2\pi/3}\frac{1}{\sqrt{1+2\cos t}}\,\dd t\displaybreak[0]\\
= {} & \frac{3}{2}\int_{-2\pi/3}^{2\pi/3}\frac{1}{\sqrt{1+2\cos t}}\,\dd t
= 3\int_0^{2\pi/3}\frac{1}{\sqrt{1+2\cos t}}\,\dd t\,.
\end{align*}
{\em Mathematica} finds 
\beqn \label{Eq:K}
  \int_0^{2\pi/3}\frac{\dd t}{\sqrt{1+2\cos t}}\,\dd t
= K\big(\sqrt{3}\big/2\big)\,,
\eeqn
where $K$ is the complete elliptic integral of the first kind \eqref{Eq:first_kind}, hence
\beqn \label{Eq:int-H_dS}
  \int_{\partial\Ol_1}H\,\dd S = 3K\big(\sqrt{3}\big/2\big)\,. 
\eeqn
A handmade proof for the relation in \eqref{Eq:K} may be found in Section \ref{S:Appendix}.\\[0.2cm]
\hspace*{0.4cm} Now we calculate the integral of mean curvature for the edges $\eps_1$, $\eps_2$ (see \eqref{Eq:M(O_1)}). Therefore, we consider $\eps_1$. The first unit normal vector $\vec{n}=\vec{n}(t)$ in a point $t\in[-2\pi/3,2\pi/3]$, $m=0$ is given by \eqref{n-vec}, the second is $\vec{n}^* = \vec{n}^*(t) =  (n_1,n_2,-n_3)^\tn{T}$. So one gets
\beq
  \alpha(t)
= \arccos\langle\vec{n}(t),\vec{n}^*(t)\rangle
= \arccos\left(\!-\frac{\cos t}{1+\cos t}\right),
\eeq
hence
\beq
  \int_{\eps_1}\alpha\,\dd t
= \int_{-2\pi/3}^{2\pi/3}\alpha(t)\,\dd t
= 2\int_0^{2\pi/3}\alpha(t)\,\dd t
= 2\int_0^{2\pi/3}\arccos\left(\!-\frac{\cos t}{1+\cos t}\right)\,\dd t\,.
\eeq 
We observe that the inverse function of the integrand is equal to the integrand, and hence the graph of the integrand symmetrical with respect to the line $f(t)=t$. As solution of
\beq
  f(t) = t = \arccos\left(\!-\frac{\cos t}{1+\cos t}\right)
\eeq
we find $t=\pi/2$, hence
\begin{align} \label{Eq:int-alpha_ds}
  \int_{\eps_1}\alpha\,\dd t
= {} & 4\int_0^{\pi/2}\left[\arccos\left(\!-\frac{\cos t}{1+\cos t}\right)-t\right]\dd t
= 4\int_0^{\pi/2}\arccos\left(\!-\frac{\cos t}{1+\cos t}\right)\dd t
- 4\int_0^{\pi/2}t\,\dd t\nonumber\\
= {} & 4\int_0^{\pi/2}\arccos\left(\!-\frac{\cos t}{1+\cos t}\right)\dd t
- \frac{\pi^2}{2}
= 4\int_0^{\pi/2}\left[\pi-\arccos\frac{\cos t}{1+\cos t}\right]\dd t
- \frac{\pi^2}{2}\nonumber\\
= {} & 4\pi\int_0^{\pi/2}\dd t - 4\int_0^{\pi/2}\arccos\frac{\cos t}{1+\cos t}\,\dd t - \frac{\pi^2}{2}\nonumber\\
= {} & \frac{3\pi^2}{2} - 4\int_0^{\pi/2}\arccos\frac{\cos t}{1+\cos t}\,\dd t\,.
\end{align}
\textit{Mathematica} and we, too, are not able to solve the last integral. It looks similar to Coxeter's integral in \cite[pp.\ 194-201]{Nahin}. The \verb+NIntegrate+-function of \textit{Mathematica} provides
\beqn \label{Eq:I}
  I:= 
\int_0^{\pi/2}\arccos\frac{\cos t}{1+\cos t}\,\dd t
\approx 1.87738105428247449505835371657\,,
\eeqn
hence
\beq
  \int_{\eps_1}\alpha\,\dd t 
\approx 7.29488238450413994801832163353\,.
\eeq
From \eqref{Eq:M(O_1)}, \eqref{Eq:int-H_dS}, \eqref{Eq:int-alpha_ds}, with \eqref{Eq:first_kind} and \eqref{Eq:I}, it follows that the integral $M$ of the mean curvature of $\Ol_r$ is given by
\beq
  M(\Ol_r)
= \left(3K\big(\sqrt{3}\big/2\big)+\frac{3\pi^2}{2}-4I\right)r
\approx 13.7644293270030696543343466299\,r\,.
\eeq
For a convex body $K$, the mean width $\bar{b}$ is given by the relation $\bar{b}(K)=M(K)/2\pi$ (see \cite[p.\ 78, Eq.\ (3.9)]{Voss}). So we we have proved the following theorem.

\begin{theorem}
The mean width of the oloid $\Ol_r$ is
\beq
  \bar{b}(\Ol_r)
= \left(\frac{3}{2\pi}\,K\big(\sqrt{3}\big/2\big)+\frac{3\pi}{4}-\frac{2}{\pi}\,I\right)r
\approx 2.19067696623158876633263049436\,r, 
\eeq
where $K$ is the complete elliptic integral of the first kind \eqref{Eq:first_kind}, and 
\beq
  I
= \int_0^{\pi/2}\arccos\frac{\cos x}{1+\cos x}\,\dd x\,. 
\eeq 
\end{theorem}

\section{Direct calculation of the mean width} \label{S:Mean_width}

Let
\beq
  P
= \{(x,y,z)\in\R^3 \colon ax+by+cz = d\} 
\eeq
be a supporting plane of $\Ol_1$ given in the Hesse normal form. So $\vec{N}=(a,b,c)^\tn{T}$ with $a,b,c\in\R$, $a^2+b^2+c^2=1$ is the normal unit vector of $P$ and $|d|$ is the distance of $P$ from the origin. $P$ intersects the plane $z=0$ in the line
\beq
  L_{xy}
= \{(x,y)\in\R^2 \colon ax+by = d\}\,, 
\eeq
and the plane $x=0$ in the line
\beq
  L_{yz}
= \{(y,z)\in\R^2 \colon by+cz = d\}\,. 
\eeq
The equation of $L_{xy}$ in Hesse normal form is
\beq
  \frac{ax}{\sqrt{a^2+b^2}} + \frac{by}{\sqrt{a^2+b^2}}
= \frac{d}{\sqrt{a^2+b^2}}\,,
\eeq
therefore, the distance $d_1$ of $L_{xy}$ from the center $(0,-1/2,0)$ of $k_A$ is
\beq
  d_1 
= \left|\frac{a}{\sqrt{a^2+b^2}}\cdot 0 + \frac{b}{\sqrt{a^2+b^2}}\cdot\left(\!-\frac{1}{2}\right)
- \frac{d}{\sqrt{a^2+b^2}}\right|
= \frac{b/2+d}{\sqrt{a^2+b^2}}
\eeq
(see e.\,g. \cite[p.\ 172]{Bosch}). Since $L_{xy}$ is tangent to $k_A$, we have
\beq
  \frac{b/2+d}{\sqrt{a^2+b^2}} = 1 \qquad\Longrightarrow\qquad
  d = \sqrt{a^2+b^2} - \frac{b}{2}\,. 
\eeq
Analogously one finds that the distance $d_2$ of $L_{yz}$ from the center $(0,1/2,0)$ of $k_B$ is
\beq
  d_2 
= \frac{b/2-d}{\sqrt{a^2+b^2}}
= 1\,,
\eeq
hence
\beq
  d = \sqrt{a^2+b^2} + \frac{b}{2}\,. 
\eeq
It follows that the distance $p$ between the support plane $P$ and the origin is  
\beq
  p = \max\left\{\sqrt{a^2+b^2}-\frac{b}{2},\sqrt{a^2+b^2}+\frac{b}{2}\right\}.
\eeq
Now we use spherical coordinates $0\le\ph\le\pi/2$ and $0\le\thet\le\pi/2$ as coordinates of the unit normal vector $\vec{N}$:
\beq
  a = \cos\ph\sin\thet\,,\quad
  b = \sin\ph\sin\thet\,,\quad
  c = \cos\thet\,.
\eeq
So we have
\begin{align*}
  \sqrt{a^2+b^2}-\frac{b}{2}
= {} & \left(1-\frac{1}{2}\,\sin\ph\right)\sin\thet\,,\\[0.2cm]
  \sqrt{a^2+b^2}+\frac{b}{2}
= {} & \frac{1}{2}\,\sin\ph\sin\thet + \sqrt{\sin^2\ph\sin^2\thet+\cos^2\thet}\,,
\end{align*}
and can write $p$ as
\beq
  p(\ph,\thet)
= \max\left\{\left(1-\frac{1}{2}\sin\ph\right)\sin\thet,
  \frac{1}{2}\sin\ph\sin\thet 
+ \sqrt{\sin^2\ph\sin^2\thet+\cos^2\thet}\right\}.
\eeq
Clearly, $p=p(\ph,\thet)$ is the support function of $\Ol_1$ in the direction $\ph,\thet$. Hence the width $b$ of $\Ol_r$ in this direction is given by
\beq
  b(\ph,\thet) = p(\ph,\thet) + p(\pi+\ph,\pi-\thet)\,.
\eeq
In order to calculate the mean width $\bar{b}$ of $\Ol_1$ we have to integrate over all directions, hence over the unit hemisphere. Let $\dd S=\dd S(\ph,\thet)=\sin\thet\,\dd\thet\,\dd\ph$ denote the surface element of the unit sphere, we have    
\begin{align*}
  \bar{b}(\Ol_1)
= {} & \frac{\int_0^\pi\int_0^\pi b(\ph,\thet)\,\dd S(\ph,\thet)}{\int_0^\pi\int_0^\pi\dd S(\ph,\thet)}
= \frac{\int_0^\pi\int_0^\pi [p(\ph,\thet)+p(\pi+\ph,\pi-\thet)]\,\dd S(\ph,\thet)}{\int_0^\pi\int_0^\pi\dd S(\ph,\thet)}\\[0.1cm]
= {} & \frac{\int_0^\pi\int_0^\pi p(\ph,\thet)\,\dd S(\ph,\thet)+\int_0^\pi\int_0^\pi p(\pi+\ph,\pi-\thet)\,\dd S(\ph,\thet)}{\int_0^\pi\int_0^\pi\dd S(\ph,\thet)}
= \frac{2\int_0^\pi\int_0^\pi p(\ph,\thet)\,\dd S(\ph,\thet)}{\int_0^\pi\int_0^\pi\dd S(\ph,\thet)}
\end{align*}
(cp.\ \cite[p.\ 78, Eq.\ (3.9)]{Voss}), where the last equal sign follows from the fact that there are two congruent portions of $\Ol_1$ in the half spaces $y\le 0$ and $y\ge 0$. Due to the symmetry of $\Ol_1$ with respect to the planes $z=0$ and $x=0$ we can restrict the spherical coordiates to the intervals $0\leq\thet\leq\pi/2$ and $0\leq\ph\leq\pi/2$, respectively, hence
\begin{align} \label{Eq:mean_width}
  \bar{b}(\Ol_1)
= {} & \frac{2\int_0^{\pi/2}\int_0^{\pi/2}p(\ph,\thet)\,\dd S(\ph,\thet)}{\int_0^{\pi/2}\int_0^{\pi/2}\dd S(\ph,\thet)}
= \frac{2\int_0^{\pi/2}\int_0^{\pi/2}p(\ph,\thet)\sin\thet\,\dd\thet\,\dd\ph}{\int_0^{\pi/2}\int_0^{\pi/2}\sin\thet\,\dd\thet\,\dd\ph}\nonumber\\[0.1cm]
= {} & \frac{4}{\pi}\int_0^{\pi/2}\int_0^{\pi/2}p(\ph,\thet)\sin\thet\,\dd\thet\,\dd\ph\nonumber\\[0.1cm]
= {} & \frac{4}{\pi}\left[\int_{\ph=0}^{\pi/6}\int_{\thet=0}^{\xi(\ph)}\left(\frac{1}{2}\sin\ph\sin\thet + \sqrt{\sin^2\ph\sin^2\thet+\cos^2\thet}\right)\sin\thet\,\dd\thet\,\dd\ph\right.\nonumber\\[0.1cm]
{} & + \int_{\ph=0}^{\pi/6}\:\int_{\thet=\xi(\ph)}^{\pi/2}\left(1-\frac{1}{2}\,\sin\ph\right)\sin^2\thet\,\dd\thet\,\dd\ph\displaybreak[0]\\[0.1cm]
{} & + \left.\int_{\ph=\pi/6}^{\pi/2}\int_{\thet=0}^{\pi/2}\left(\frac{1}{2}\sin\ph\sin\thet + \sqrt{\sin^2\ph\sin^2\thet+\cos^2\thet}\right)\sin\thet\,\dd\thet\,\dd\ph\right],\nonumber
\end{align}
where
\beq
  \xi(\ph)
= \arccos\frac{\sqrt{-1+2\sin\ph}}{\sqrt{-2+2\sin\ph}}
\eeq
is the solution of the equation
\beq
  \left(1-\frac{1}{2}\sin\ph\right)\sin\thet
= \frac{1}{2}\sin\ph\sin\thet + \sqrt{\sin^2\ph\sin^2\thet+\cos^2\thet}
\eeq
for $\thet=\xi(\ph)$. Numerical integration of \eqref{Eq:mean_width} with {\em Mathematica} gives
\beq
  \bar{b}(\Ol_1) \approx 2.19067696623\,. 
\eeq


\section{The parallel body} \label{S:Parallel_body}

For a convex body $K\subset\R^n$ and $\rh>0$, the set (Minkowski sum)
\beq
  K+B_\rh^n = \{x\in\R^n\colon d(x,K)\le\rh\}
\eeq
is the {\em parallel body} of $K$ at distance $\rh$, where $B_\rh^n$ is the $n$-ball of radius $\rh$,
\beq
  B_\rh^n=\{x\in\R^n\colon\lVert x\rVert\leq\rh\}\,,
\eeq
and $d(x,K)$ is the distance between the point $x$ and $K$. The volume of the parallel body is given by the Steiner formula
\beqn \label{Eq:vol_parallel_body}
  V_n(K+B_\rh^n)
= \sum_{j=0}^n\rh^{n-j}\kappa_{n-j}V_j(K)\,,
\eeqn 
where
\beqn \label{Eq:kappa}
  \kappa_k
= \frac{\pi^{k/2}}{\varGamma(1+k/2)}
\eeqn
is the volume of the $k$-dimensional unit ball $B_1^k$, and $V_0,\ldots,V_{n-1}$ are the intrinsic volumes of $K$. \cite[p.\ 2,\ pp.\ 12-13, p.\ 600]{SchneiderWeil3} \\[0.2cm]
\hspace*{0.4cm} Using the relations in \cite[p.\ 301]{SchneiderWeil2}, where $\chi$ denotes the Euler characteristic, the intrinsic volumes of $\Ol_r$ are
\beqn \label{Eq:intrinsic_volumes}
\left.
\begin{aligned}
  V_0(\Ol_r) 
= {} & \chi = 1\,,\quad
  V_1(\Ol_r) 
= 2\bar{b}(\Ol_r)
= \frac{M(\Ol_r)}{\pi}
= \left[\frac{3}{\pi}\,K\big(\sqrt{3}\big/2\big)+\frac{3\pi}{2}-\frac{4}{\pi}\,I\right]r\,,\\[0.1cm]
  V_2(\Ol_r)
= {} & \frac{1}{2}\,S(\Ol_r)
= 2\pi r^2\,,\quad
  V_3(\Ol_r)
= V(\Ol_r)
= \frac{2}{3}\left[2E\big(\sqrt{3}\big/2\big)+K\big(\sqrt{3}\big/2\big)\right]r^3 
\end{aligned}
\;\;\right\}
\eeqn
with $E$ (see \eqref{Eq:second_kind}), $K$ (see \eqref{Eq:first_kind}), and $I$ (see \eqref{Eq:I}). Let $\Ol_{r,\,\rh}$ denote the parallel body of $\Ol_r$ at distance $\rh$. Due to \eqref{Eq:vol_parallel_body}, its volume is
\begin{align*}
  V(\Ol_{r,\,\rh}) = V_3(\Ol_{r,\,\rh})
= {} & \kappa_0V_3(\Ol_r) + \kappa_1V_2(\Ol_r)\rh + \kappa_2V_1(\Ol_r)\rh^2 + \kappa_3V_0(\Ol_r)\rh^3\\
= {} & V_3(\Ol_r) + 2V_2(\Ol_r)\rh + \pi V_1(\Ol_r)\rh^2 + \frac{4\pi}{3}V_0(\Ol_r)\rh^3\\
= {} & V(\Ol_r) + S(\Ol_r)\rh + M(\Ol_r)\rh^2 + \frac{4\pi}{3}\rh^3\\
= {} & V(\Ol_1)r^3 + S(\Ol_1)r^2\rh + M(\Ol_1)r\rh^2 + \frac{4\pi}{3}\rh^3\,.
\end{align*}
Applying \cite[p.\ 82, (3.17)]{Voss} allows to calculate the surface area $S$ of $\Ol_{r,\,\rh}$:
\beq
  S(\Ol_{r,\,\rh})
= S(\Ol_r) + 2M(\Ol_r)\rh + 4\pi\rh^2
= S(\Ol_1)r^2 + 2M(\Ol_1)r\rh + 4\pi\rh^2\,.
\eeq
Clearly, the mean width of $\Ol_{r,\,\rh}$ is equal to $\bar{b}(\Ol_r)+2\rh$, hence
\beq
  M(\Ol_{r,\,\rh})
= 2\pi\left[\bar{b}(\Ol_r)+2\rh\right]
= 2\pi\bar{b}(\Ol_r)+4\pi\rh
= M(\Ol_r)+4\pi\rh
= M(\Ol_1)r+4\pi\rh
\eeq
(see also \cite[p.\ 82, (3.17)]{Voss}). The results of the following theorem follow immediately.
\begin{theorem}
The integral $M$ of mean curvature, the surface area $S$ and the volume $V$ of the parallel body $\Ol_{r,\,\rh}$ are given by
\begin{align*}
  M(\Ol_{r,\,\rh})
= {} & M(\Ol_1)\,r 
+ 4\pi\rh\,,\qquad
  S(\Ol_{r,\,\rh})
= 4\pi r^2
+ 2M(\Ol_1)\,r\rh
+ 4\pi\rh^2\,,\\[0.2cm]
  V(\Ol_{r,\,\rh})
= {} & \frac{2}{3}\left[2E\big(\sqrt{3}\big/2\big)+K\big(\sqrt{3}\big/2\big)\right]r^3
+ 4\pi r^2\rh 
+ M(\Ol_1)\,r\rh^2
+ \frac{4\pi}{3}\,\rh^3
\end{align*}
with
\beq
  M(\Ol_1)
= 3K\big(\sqrt{3}\big/2\big) + \frac{3\pi^2}{2} 
- 4\int_0^{\pi/2}\arccos\frac{\cos x}{1+\cos x}\,\dd x\,.
\eeq
\end{theorem}

\section{Intersections with an oloid} \label{S:Intersections}

Now, we apply our results and the principal kinematic formula to derive some expectations for the intersections of the oloid $\Ol_r$ and the three-dimensional ball $B_r:=B_r^3$ of radius $r$, and of two oloids $\Ol_r$.\\[0.2cm]
\hspace*{0.4cm} The principal kinematic formula (see \cite[p.\ 301]{SchneiderWeil2}) for a fixed convex body $K$ and a moving convex body $M$ is for $j\in\{0,\ldots,n\}$ given by 
\beqn \label{Eq:kinematic_formula}
  I_j(K,M)
:= \int_{SO_n}\int_{\R^n}V_j(K\cap(\thet M+\vec{x}))\,\dd\lambda(\vec{x})\,\dd\nu(\thet)
= \sum_{k=j}^n\alpha_{njk}V_k(K)V_{n+j-k}(M)
\eeqn
with the notation
\begin{center}
\begin{tabular}{ll}
$SO_n$    & group of proper (orientation-preserving) rotations \cite[p.\ 13]{SchneiderWeil3},\\
$\thet$   & proper rotation, $\thet\in SO_n$,\\
$\vec{x}$ & translation vector,\\
$\lambda$ & Lebesgue measure on $\R^n$,\\
$\nu$     & unique Haar measure on $SO_n$ with $\nu(SO_n)=1$ \cite[p.\ 584]{SchneiderWeil3},
\end{tabular}
\end{center}  
and
\beq
  \alpha_{njk}
= \frac{k!\,\kappa_k(n+j-k)!\,\kappa_{n+j-k}}{j!\,\kappa_j\,n!\,\kappa_n}\,,\quad
  \alpha_{njk}
= \alpha_{nj(n+j-k)}\,,\quad
  \alpha_{njj}
= \alpha_{njn}
= 1\,, 
\eeq
where $\kappa_k$ is the volume of the unit $k$-ball (see \eqref{Eq:kappa}). Since the intersection of two convex sets is a convex set, we have
\begin{align} \label{Eq:I_0(K,M)}
  I_0(K,M)
= {} & \int_{SO_n}\int_{\R^n}\chi(K\cap(\thet M+\vec{x}))\,\dd\lambda(\vec{x})\,\dd\nu(\thet)\nonumber\\[0.1cm]
= {} & \int_{SO_n}\int_{\R^n}\mathbbm{1}_{K\,\cap\,(\thet M\,+\,\vec{x})\,\ne\,\emptyset}\:\dd\lambda(\vec{x})\,\dd\nu(\thet)\,,
\end{align}
where $\mathbbm{1}_B$ is the indicator function of the event $B$. So we see that $I_0(K,M)$ is the measure of the set of rigid motions bringing $M$ into a hitting position with $K$ (see \cite[p.~175]{SchneiderWeil3}, \cite[p.\ 262, p.\ 267]{Santalo}).\\[0.2cm] 
\hspace*{0.4cm} For $n=3$, \eqref{Eq:kinematic_formula} gives
\beqn \label{Eq:I_j(K,M)}
\left.
\begin{aligned}
  I_0(K,M)
= {} & V_0(K)\,V_3(M) + \frac{1}{2}\,V_1(K)\,V_2(M) 
  + \frac{1}{2}\,V_2(K)\,V_1(M) + V_3(K)\,V_0(M)\,,\\[0.1cm]		
  I_1(K,M)
= {} & V_1(K)\,V_3(M) + \frac{\pi}{4}\,V_2(K)\,V_2(M) 
  + V_3(K)\,V_1(M)\,,\\[0.1cm]
  I_2(K,M)
= {} & V_2(K)\,V_3(M) + V_3(K)\,V_2(M)\,,\\[0.1cm]
  I_3(K,M)
= {} & V_3(K)\,V_3(M)\,.		
\end{aligned}
\;\;\right\}
\eeqn
From \eqref{Eq:I_0(K,M)} it follows that
\begin{align} \label{Eq:E[V]}
  \E\left[V(K\cap M)\right]
= {} &  \E\left[V_3(K\cap M)\right]
= \frac{I_3(K,M)}{I_0(K,M)}\,.
\end{align}
is the expected volume of $K\cap M$. Analogously, we get the expected mean width and the expected surface area:   
\begin{align}
  \E\big[\bar{b}\big(K\cap M\big)\big]
= {} & \frac{1}{2}\,\E\big[V_1(K\cap M)\big]
= \frac{I_1(K,M)}{2I_0(K,M)}\,,\label{Eq:E[b]}\\[0.1cm]
  \E\big[S\big(K\cap M\big)\big]
= {} & 2\,\E\big[V_2(K\cap M)\big]
= \frac{2I_2(K,M)}{I_0(K,M)}\,.\label{Eq:E[S]}
\end{align}
Clearly, it is possible to reverse the roles of the fixed body and the moving body.

\begin{example}
As an example we calculate the expected values \eqref{Eq:E[b]}, \eqref{Eq:E[S]} and \eqref{Eq:E[V]} for $K=\Ol_r$ and $M=B_r$, or, equivalently, for $K=B_r$ and $M=\Ol_r$. For the ball $B_r$ one easily gets
\beqn \label{Eq:intrinsic_volumes_ball}
\left.
\begin{aligned}
  V_0(B_r)
= {} & \chi(B_r) 
= 1\,,\quad
  V_1(B_r)
= 2\bar{b}(B_r)
= 4r\,,\\
  V_2(B_r)
= {} & \frac{1}{2}\,S(B_r)
= 2\pi r^2\,,\quad
  V_3(B_r)
= V(B_r)
= \frac{4\pi r^3}{3}\,.  
\end{aligned}
\;\;\right\}
\eeqn
Note that these terms also follow from the general formula \beq
  V_k(B_r)
= V_k\big(B_r^3\big)
= V_k(B_1^3)\,r^k
= {{3}\choose{k}}\,\frac{\kappa_3}{\kappa_{3-k}}\,r^k
\eeq
\cite[p.\ 300]{SchneiderWeil2}, where $\kappa_k$ is the volume of the unit $k$-ball (see \eqref{Eq:kappa}). Plugging \eqref{Eq:intrinsic_volumes} and \eqref{Eq:intrinsic_volumes_ball} in \eqref{Eq:I_j(K,M)} gives
\begin{align*}
  I_0(\Ol_r,B_r)
= {} & \frac{r^3}{6}\left[9\pi^2 + 32\pi + 8E\big(\sqrt{3}\big/2\big)
  + 22K\big(\sqrt{3}\big/2\big)-24I\right],\\[0.1cm]
  I_1(\Ol_r,B_r)
= {} & \frac{r^4}{3}\left[3\pi^3 + 6\pi^2 + 16E\big(\sqrt{3}\big/2\big)
  + 20K\big(\sqrt{3}\big/2\big)-16I\right],\\[0.1cm]
  I_2(\Ol_r,B_r)
= {} & \frac{4\pi r^5}{3}\left[2\pi + 2E\big(\sqrt{3}\big/2\big)
  + K\big(\sqrt{3}\big/2\big)\right],\\[0.1cm]
  I_3(\Ol_r,B_r)
= {} & \frac{8\pi r^6}{9}\left[2E\big(\sqrt{3}\big/2\big)
  + K\big(\sqrt{3}\big/2\big)\right],
\end{align*}
and
\begin{align*}
  \E\big[\bar{b}(\Ol_r\cap B_r)\big]
= {} & \frac{I_1(\Ol_r,B_r)}{2I_0(\Ol_r,B_r)}
= \frac{3\pi^3 + 6\pi^2 + 16E\big(\sqrt{3}\big/2\big)
  + 20K\big(\sqrt{3}\big/2\big)-16I}
  {9\pi^2 + 32\pi + 8E\big(\sqrt{3}\big/2\big)
  + 22K\big(\sqrt{3}\big/2\big)-24I}\,r\,,\\[0.1cm]
  \E\big[S(\Ol_r\cap B_r)\big]
= {} & \frac{2I_2(\Ol_r,B_r)}{I_0(\Ol_r,B_r)}
= \frac{16\pi\left[2\pi + 2E\big(\sqrt{3}\big/2\big)
  + K\big(\sqrt{3}\big/2\big)\right]}
  {9\pi^2 + 32\pi + 8E\big(\sqrt{3}\big/2\big)
  + 22K\big(\sqrt{3}\big/2\big)-24I}\,r^2\,,\\[0.1cm]
  \E\big[V(\Ol_r\cap B_r)\big]
= {} & \frac{I_3(\Ol_r,B_r)}{I_0(\Ol_r,B_r)}
= \frac{16\pi\left[2E\big(\sqrt{3}\big/2\big)
  + K\big(\sqrt{3}\big/2\big)\right]}
  {3\left[9\pi^2 + 32\pi + 8E\big(\sqrt{3}\big/2\big)
  + 22K\big(\sqrt{3}\big/2\big)-24I\right]}\,r^3\,.
\end{align*}
\end{example}

\begin{example}
In the case $K=\Ol_r=M$, we have
\begin{align*}
  I_0(\Ol_r,\Ol_r)
= {} & \frac{r^3}{3}\left[9\pi^2 + 8E\big(\sqrt{3}\big/2\big)
  + 22K\big(\sqrt{3}\big/2\big)-24I\right],\\[0.1cm]
  I_1(\Ol_r,\Ol_r)
= {} & \frac{r^4}{3\pi}\left[3\pi^4 + 2\left(2E\big(\sqrt{3}\big/2\big)
  + K\big(\sqrt{3}\big/2\big)\right)
  \left(3\pi^2+6K\big(\sqrt{3}\big/2\big)-8I\right)\right],\\[0.1cm]
  I_2(\Ol_r,\Ol_r)
= {} & \frac{8\pi r^5}{3}\left[2E\big(\sqrt{3}\big/2\big)
  + K\big(\sqrt{3}\big/2\big)\right],\\[0.1cm]
  I_3(\Ol_r,\Ol_r)
= {} & \frac{4r^6}{9}\left[2E\big(\sqrt{3}\big/2\big)
  + K\big(\sqrt{3}\big/2\big)\right]^2,
\end{align*}
hence
\begin{align*}
  \E\big[\bar{b}(\Ol_r\cap\Ol_r)\big]
= {} &  \frac{I_1(\Ol_r,\Ol_r)}{2I_0(\Ol_r,\Ol_r)} 
=  \frac{3\pi^4 + 2\left[2E\big(\sqrt{3}\big/2\big)
  + K\big(\sqrt{3}\big/2\big)\right]
  \left[3\pi^2+6K\big(\sqrt{3}\big/2\big)-8I\right]}
  {2\pi\left[9\pi^2 + 8E\big(\sqrt{3}\big/2\big)
  + 22K\big(\sqrt{3}\big/2\big)-24I\right]}\,r\,,\\[0.1cm]
  \E\big[S(\Ol_r\cap\Ol_r)\big]
= {} & \frac{2I_2(\Ol_r,\Ol_r)}{I_0(\Ol_r,\Ol_r)}
= \frac{16\pi\left[2E\big(\sqrt{3}\big/2\big)
  + K\big(\sqrt{3}\big/2\big)\right]}
  {9\pi^2 + 8E\big(\sqrt{3}\big/2\big)
  + 22K\big(\sqrt{3}\big/2\big)-24I}\,r^2\,,\\[0.1cm]
  \E\big[V(\Ol_r\cap\Ol_r)\big]
= {} & \frac{I_3(\Ol_r,\Ol_r)}{I_0(\Ol_r,\Ol_r)}
= \frac{4\left[2E\big(\sqrt{3}\big/2\big)
  + K\big(\sqrt{3}\big/2\big)\right]^2}
  {3\left[9\pi^2 + 8E\big(\sqrt{3}\big/2\big)
  + 22K\big(\sqrt{3}\big/2\big)-24I\right]}\,r^3\,.
\end{align*}
\end{example}

\noindent
The following table shows numerical approximations for the expectations of the intersections.

\begin{center}
\renewcommand{\arraystretch}{1.25}
\begin{tabular}{|c|c||c|c|c|} \hline
$K$ & $M$ & $\E\big[\bar{b}(K\cap M)\big]\big/r$ & 
  $\E\big[S(K\cap M)\big]\big/r^2$ & $\E\big[V(K\cap M)\big]\big/r^3$\\ \hline\hline
 $B_r$  &  $B_r$  & 0.9626377063 & 3.141592654 & 0.5235987756\\
$\Ol_r$ &  $B_r$  & 0.9169621588 & 2.710463736 & 0.3808512243\\
$\Ol_r$ & $\Ol_r$ & 0.8585694641 & 2.280916270 & 0.2770215506\\ \hline
\end{tabular}
\renewcommand{\arraystretch}{1}
\end{center}
\vspace{0.1cm}


\section{Appendix} \label{S:Appendix}

Now we are going to show that the integral
\beq
  J
:= \int_0^{2\pi/3}\frac{\dd t}{\sqrt{1+2\cos t}}
\eeq
is equal to $K\big(\sqrt{3}\big/2\big)$ (see \eqref{Eq:K}) without the use of {\em Mathematica}. With
\beq
  1+2\cos t
= 1+2\left(1-2\sin^2\frac{t}{2}\right)
= 3-4\sin^2\frac{t}{2}
= 3\left[1-\left(\frac{2}{\sqrt{3}}\right)^{\!2}\sin^2\frac{t}{2}\right]
\eeq
we have
\beq
  J
= \int_0^{2\pi/3}\frac{\dd t}{\sqrt{3}\:\sqrt{1-\big(2\big/\sqrt{3}\,\big)^2\sin^2(t/2)}}\,.
\eeq
Now, following the argumentation in \cite{Weisstein}, we put
\beq
  \csc\frac{t_0}{2} = \frac{2}{\sqrt{3}}
  \qquad\Longrightarrow\qquad
  t_0 = \frac{2\pi}{3}\,.
\eeq
This gives
\beq
  J
= \frac{1}{2}\int_0^{t_0}\frac{1}{\sin(t_0/2)\,\sqrt{1-\csc^2(t_0/2)\sin^2(t/2)}}\,\dd t\,.
\eeq
Now let
\beq
  \sin(t/2)
= \sin(t_0/2)\,\sin\ph\,,
\eeq
so the angle $t$ is transformed to
\beq
  \ph 
= \arcsin\frac{\sin(t/2)}{\sin(t_0/2)}\,,
\eeq
hence
\beq
  t = 0 \;\;\Longrightarrow\;\; \ph=0\,,\qquad
  t = t_0 \;\;\Longrightarrow\;\; \ph=\frac{\pi}{2}\,.
\eeq
Taking the differentials gives
\beq
  \frac{1}{2}\cos\left(\frac{t}{2}\right)\dd t
= \sin\left(\frac{t_0}{2}\right)\cos\ph\,\dd\ph\,,
\eeq
or
\beq
  \frac{1}{2}\:\sqrt{1-\sin^2\left(\frac{t}{2}\right)}\;\dd t
= \sin\left(\frac{t_0}{2}\right)\cos\ph\,\dd\ph\,,
\eeq
hence
\beq
  \frac{1}{2}\:\sqrt{1-\sin^2\left(\frac{t_0}{2}\right)\sin^2\ph}\;\dd t
= \sin\left(\frac{t_0}{2}\right)\cos\ph\,\dd\ph\,.
\eeq
Plugging this in gives
\begin{align*}
  J
= {} & \frac{1}{2}\int_0^{\pi/2}
  \frac{1}{\sin(t_0/2)\,\sqrt{1-\sin^2\ph}}\,
  \frac{\sin(t_0/2)\cos\ph\,\dd\ph}
       {(1/2)\,\sqrt{1-\sin^2(t_0/2)\sin^2\ph}}\\
= {} & \int_0^{\pi/2}
  \frac{\dd\ph}{\sqrt{1-\sin^2(t_0/2)\sin^2\ph}}
= \int_0^{\pi/2}
  \frac{\dd\ph}{\sqrt{1-\sin^2(\pi/3)\sin^2\ph}}\\
= {} & \int_0^{\pi/2}
  \frac{\dd\ph}{\sqrt{1-\big(\sqrt{3}\big/2\big)^2\sin^2\ph}}
= K\big(\sqrt{3}\big/2\big)\,.  
\end{align*}
 


\bigskip

\noindent
Author's address:\\[0.2cm]
Uwe B\"asel\\ 
HTWK Leipzig, University of Applied Sciences,\\
Faculty of Mechanical and Energy Engineering,\\
PF 30 11 66, 04251 Leipzig, Germany\\
e-mail: uwe.baesel@htwk-leipzig.de

\end{document}